\documentclass[11pt]{amsproc}
\usepackage{a4}
\usepackage[dvips]{graphicx}
\usepackage{amsmath}
\usepackage{amsthm}
\usepackage{amssymb}
\usepackage{amsfonts}
\usepackage{mathrsfs}
\usepackage[all]{xy}

\theoremstyle{plain}
\newtheorem{theorem}[subsubsection]{Theorem}

\newtheorem{lemma}[subsubsection]{Lemma}

\newtheorem{proposition}[subsubsection]{Proposition}

\theoremstyle{definition}
\newtheorem{definition}[subsubsection]{Definition}
\newtheorem{definition-subsection}[subsection]{Definition}

\theoremstyle{remark}
\newtheorem{aside}[subsubsection]{Aside}
\newtheorem{remark}[subsubsection]{Remark}
\newtheorem{example}[subsubsection]{Example}

\newcommand{\Oc}{\mathcal{O}}

\newcommand{\Mcal}{\mathcal{M}}
\newcommand{\Sc}{\mathcal{S}}
\newcommand{\Pc}{\mathcal{P}}
\newcommand{\Ec}{\mathcal{E}}

\newcommand{\Rc}{\mathcal{R}}
\newcommand{\Tc}{\mathcal{T}}

\newcommand{\Ac}{\mathcal{A}}

\newcommand{\Db}{\mathbb{D}}
\newcommand{\A}{\mathbb{A}}
\newcommand{\C}{\mathbb{C}}
\newcommand{\R}{\mathbb{R}}
\newcommand{\Q}{{\mathbb{Q}}}

\newcommand{\Z}{\mathbb{Z}}

\newcommand{\prelilacs}{\textrm{(\textbf{pre-Lilacs})}}

\newcommand{\G}{\mathbb{G}}
\newcommand{\Sb}{\mathbb{S}}
\newcommand{\Hb}{\mathbb{H}}
\newcommand{\Pb}{\mathbb{P}}

\newcommand{\Ker}{\mathrm{Ker}}

\newcommand{\Hom}{\mathrm{Hom}}
\newcommand{\Pic}{\mathrm{Pic}}

\newcommand{\im}{\mathrm{im}}
\newcommand{\card}{\mathrm{card}}
\newcommand{\End}{\mathrm{End}}
\newcommand{\diag}{\mathrm{diag}}

\newcommand{\id}{\mathrm{id}}
\newcommand{\ch}{\mathrm{ch}}

\newcommand{\limproj}{\underset{\leftarrow}{\mathrm{lim}}\,}

\newcommand{\PGL}{\mathrm{PGL}}

\newcommand{\Irrat}{\mathrm{Irrat}}

\newcommand{\HC}{\mathrm{HC}}
\newcommand{\GL}{\mathrm{GL}}

\newcommand{\nct}{\mathrm{(\mathbf{Nct})}}

\newcommand{\GSp}{\mathrm{GSp}}
\newcommand{\MT}{\mathrm{MT}}
\newcommand{\BMT}{\mathrm{BMT}}

\newcommand{\Sh}{\mathrm{Sh}}

\newcommand{\Res}{\mathrm{Res}}

\title[Three examples of NC boundaries of Shimura varieties]{Three 
examples of non-commutative boundaries of Shimura varieties}
\author{Fr\'ed\'eric Paugam}
\email{fpaugam@math.univ-rennes1.fr}
\date{}

\begin{document}
\maketitle

\begin{abstract}
Our modest aims in writing this paper were twofold:
we first wanted to understand the linear algebra and algebraic group
theoretic background of Manin's real multiplication program proposed in
\cite{Manin3}.
Secondly, we wanted to find nice higher dimensional analogs
of the non-commutative modular curve studied by Manin and Marcolli in
\cite{Manin-Marcolli2}.
These higher dimensional objects, that we call \emph{irrational} or 
\emph{non-commutative boundaries of Shimura varieties},
are double cosets spaces of the form
$\Gamma\backslash G(\R)/P(K)$,
where $G$ is a (connected) reductive $\Q$-algebraic group,
$P(K)=M(K)AN\subset G(\R)$ is
a real parabolic subgroup
corresponding to a rational parabolic subgroup
$P\subset G$, and $\Gamma\subset G(\Q)$ is an arithmetic subgroup.
Along the way, it also seemed clear that the spaces
$$\Gamma\backslash G(\R)/M(K)A$$
are of great interest, and sometimes more convenient to study.
We study in this document three examples of these general
spaces. These spaces describe degenerations of complex structures
on tori in (multi)foliations.
\end{abstract}


\section*{Introduction}
The non-commutative modular curve is the chaotic space
$$\GL_2(\Z)\backslash\Pb^1(\R).$$
Its first appearance in the non-commutative geometric world arose in the work
\cite{Connes4} of Connes on differential geometry of non-commutative tori.
Since the action of $\GL_2(\Z)$ by homographies on $\Pb^1(\R)$ is very
chaotic, the classical quotient space, whose algebra is the one of continuous
complex functions on $\Pb^1(\R)$ that are invariant by $\GL_2(\Z)$, is 
topologically
identified with a point. The philosophy of Connes' non-commutative geometry
and the related philosophy of topoi (as explained by Cartier in
\cite{Cartier}) tell us that this
quotient space is much more than that.
Connes showed that the crossed product $C^*$-algebra
$C^*(\Pb^1(\R))\rtimes \GL_2(\Z)$
is a good analog of the algebra of continuous functions
for such a chaotic space because
it is possible to calculate from it nice cohomological invariants, as $K$-theory
and cyclic cohomology (from its $C^\infty$ version), that have a real geometric
meaning. He showed, indeed, that
cyclic cohomology is a good analog of De Rham cohomology of manifolds for such
chaotic spaces. This cohomology theory is a very profound tool that permitted
Connes and Moscovici,
to cite one example among many others, to prove a local index
formula for foliations, which is an analog of the Riemann-Roch theorem for
the spaces of leaves \cite{Connes5}. This result was not accessible without
non-commutative geometric intuition. Bost-Connes and Connes have also shown 
\cite{Bost-Connes},
\cite{Connes6},  that non-commutative geometry can be very useful for 
arithmetic questions. We refer the reader
to Marcolli \cite{Mat4} for a nice survey on non-commutative arithmetic geometry.

The non-commutative modular curve also appeared some years ago
in Connes, Douglas and Schwarz paper \cite{Connes1}
as moduli space for physical backgrounds for
compactifications of string theory .
It was already shown in \cite{Seiberg-Witten} by Seiberg and Witten that 
non-commutative
spaces can be good backgrounds for open string theory. The deformation 
quantization story
of Kontsevich and Soibelman also gives related results.
The author, being more informed of arithmetic geometry, will not discuss 
these
physical motivations that are very important for future developments of 
non-commutative
moduli spaces.

The arithmetic viewpoint of the non-commutative modular curve appeared first 
in
the work of Manin and Marcolli \cite{Manin-Marcolli2}, \cite{Manin-Marcolli}
and in Manin's real multiplication program
\cite{Manin3}. These works were for us the main inspirations for writing this
paper.

The basic idea in our work is the following: the classical theory of Shimura 
varieties
was completely rewritten in terms of Hodge structures by Deligne in \cite{De4}
and it proved to be very helpful for arithmetics, for example in the theory
of absolute Hodge motives made in \cite{De1} and for the construction of
canonical models for all Shimura varieties, made by Milne in \cite{Milne5}.
One of the interests of this construction is to translate fine information
about moduli spaces in terms of algebraic groups morphisms. The author
wanted to know if it was possible to make such a translation for non-commutative
boundaries of Shimura varieties.

As higher dimensional analogs of the boundary $\Pb^1(\R)$
of the double half plane $\Hb^\pm:=\Pb^1(\C)-\Pb^1(\R)$, we
chose to use the components
of Satake's compactifications of symmetric spaces.
These components
can be written as quotients $G(\R)/P(K)$,
where $G$ is a (connected) reductive $\Q$-algebraic group and
$P(K)=M(K)AN\subset G(\R)$ is
a real parabolic subgroup\footnote{See
the book \cite{BoJi} for the definition and decomposition
of these real parabolic subgroups.} corresponding to a rational
parabolic subgroup $P\subset G$. The higher dimensional analogs
of the non-commutative modular curve $\GL_2(\Z)\backslash\Pb^1(\R)$,
which we call \emph{irrational} or 
\emph{non-commutative boundaries of Shimura varieties}, are given
by double coset spaces
$$\Gamma\backslash G(\R)/P(K),$$
where $\Gamma\subset G(\Q)$ is an arithmetic subgroup.
Along the way, it seemed also clear that the spaces
$$\Gamma\backslash G(\R)/M(K)A,$$
which we call \emph{irrational} or
\emph{non-commutative shores of Shimura varieties},
are of great interest, and sometimes more convenient to study from
the algebraic group theoretical viewpoint.

The plan of this document is the following.
In the first part, we study special geodesics on the modular curve that are
good analogs of elliptic curves with complex multiplication. We explain
how to relate the counting of these geodesics to number theoretical 
considerations.
In the second part, we construct the moduli space and universal family
of non-commutative tori in a way analogous to the construction of the
moduli space of elliptic curves and its universal family.
In the third part, we give two higher dimensional examples of non-commutative 
shores of Shimura varieties that parametrize some
degenerations of complex structures on tori in multi-foliations.

\section*{Acknowledgments}
I thank Matilde Marcolli for introducing me to non-commutative geometry and
sharing with me some of her deep insights into
non-commutative arithmetic geometry.

I thank the following people for interesting discussions and/or valuable help
in the preparation of this paper: Y. Andr\'e,
D. Blottiere, J.-B. Bost, B. Calmes, S. Cantat,
A. Connes,
F. Dal'bo, E. Ha, E. Ghys, Y. Guivarch', K. K\"unnemann,
Y. Manin, R. Noot, N. Ramachandran,
A. Thuillier, D. Zagier.
I thank the referee for useful comments and nice suggestions.

I thank the following institutions and workshops for financial
and technical support:
Bonn's Max Planck Institut f\"ur Mathematik,
Bonn's worshops ``Number theory and non-commutative geometry I and II'',
the Les Houches school ``Fronti\`eres entre th\'eorie des nombres,
physique et g\'eom\'etrie'',
Regensburg's ``Oberseminar Nichtkommutativ Geometrie'',
Regensburg's mathematics laboratory,
Rennes's mathematics laboratory IRMAR,
Rennes's workshop ``Jeunes chercheurs en g\'eom\'etrie'',
the European networks ``K-theory and Algebraic Groups''
and ``Arithmetic Algebraic Geometry''.

\section*{Notations}
\begin{definition-subsection}
A pair $(P,X)$ consisting of a $\Q$-algebraic group $P$
and a left $P(\R)$-space $X$ is called a \emph{pre-Shimura datum}\footnote{Most
of the pre-Shimura data that will appear in this
document will be constructed using conjugacy classes of morphisms. However,
the target group will not always be $P$ and some useful morphisms
between them will not be morphisms of conjugacy classes but
just equivariant morphisms. This is why we use such a weak definition.}.
A
\emph{morphism of pre-Shimura data}
$(P_1,X_1)\to (P_2,X_2)$ is a pair $(\phi,\psi)$ consisting
of a morphism $\phi:P_1\to P_2$ of groups and a $P_1(\R)$-equivariant map
$\psi:X_1\to X_2$.
If $(P,X)$ is a pre-Shimura datum
and $K\subset P(\A_f)$ is a compact open subgroup, the set
$$\Sh_K(P,X):=P(\Q)\backslash (X\times P(\A_f)/K)$$
is called the \emph{pre-Shimura set of level $K$} for $(P,X)$,
and the set
$$\Sh(P,X):=\limproj_K \Sh_K(P,X),$$
where $K$ runs over all compact open subgroups in $P(\A_f)$,
is called the \emph{pre-Shimura tower} associated to $(P,X)$.
\end{definition-subsection}

We warn the reader that even if $X$ has a nontrivial
$C^\infty$-structure, the corresponding Shimura space, viewed as
a quotient topological space,
can be very degenerate (even trivial).
In such cases, it may be more interesting to study the crossed
product algebra
$$C^\infty(X\times P(\A_f)/K)\rtimes P(\Q).$$

These are essentially the kind of non-commutative spaces we will
consider in this paper.

\section{Special geodesics on the modular curve}
The modular curve is one of the basic examples of a Shimura variety.
Its interpretation as a moduli space of Hodge structures (or simpler
of complex structures) allows one to nicely define \emph{special points}
of this curve as those whose Mumford-Tate group is a torus.

We propose an analogous construction for the space
$\GL_2(\Z)\backslash\GL_2(\R)/D(\R)$ of leaves of the geodesic foliation
(where $D(\R)$ is the subgroup of diagonal matrices in $\GL_2(\R)$).
This naive translation gives another point of view of the strong analogy
between the modular curve and the space of leaves of the geodesic flow,
which was already known to Gauss. Being closer to the modern (i.e., adelic)
point of view of Shimura varieties increases its chances to be generalized
to boundaries of higher dimensional Shimura varieties.

\subsection{Mumford-Tate groups of elliptic curves}
Let $E$ be a complex elliptic curve. Let $M=H^1(E,\Q)$ be its first singular
homology group. We have the Hodge decomposition
$$M_\C=H^1(E,\C)\cong H^0(E,\Omega^1_E)\oplus H^1(E,\Oc_E).$$
Let $\G_{m,\C}$ be the multiplicative group of $\C$ as an algebraic group.
We let $(x,y)\in (\C^*)^2$ act on $M_\C$ by multiplication by $x$ on
$H^0(E,\Omega^1_E)$ and $y$ on $H^1(E,\Oc_E)$.
This defines a natural morphism
$h:\G_{m,\C}^2\to \GL(M_\C)$. The \emph{Mumford-Tate group} of $E$ is the
smallest $\Q$-algebraic subgroup of $\GL(M)$ that contains the image of
$h$ over $\C$.
The following proposition serves us as a guide to define Mumford-Tate groups
of geodesics.

\begin{proposition}
Let $E$ be an elliptic curve over $\C$. The Mumford-Tate group of $E$ is
either $\GL_2$, or $\Res_{K/\Q}\G_m$ (i.e., the group $K^\times$ as a 
$\Q$-algebraic group) with $K/\Q$ an imaginary quadratic field.
In this case, we say that the curve is \emph{special} or
\emph{with complex multiplication}.
\end{proposition}

\subsection{Geodesics and analogs of Hodge structures}
\label{analogs-hodge-structures}
We recall that the space of geodesics on the modular curve can be written as
a double coset space
$Y:=\GL_2(\Z)\backslash \GL_2(\R)/D(\R)$
with $D(\R)$ the subgroup of diagonal matrices in $\GL_2(\R)$.
We denote by $\G_{m,\R}$ the multiplicative group $\R^*$ viewed as an
algebraic group.
Let $h_0:\G_{m,\R}^2\to \GL_{2,\R}$ be the morphism of algebraic
groups
that sends a pair $(x,y)\in(\R^*)^2$ to the diagonal matrix
$\diag(x,y):=
\left(\begin{smallmatrix}
x & 0\\
0 & y
\end{smallmatrix}\right)$.
We remark that $\GL_2(\R)$ acts on the space $\Hom(\G_{m,\R}^2,\GL_{2,\R})$ by
conjugation. Let
$$X:=\GL_2(\R)\cdot h_0:=\{gh_0g^{-1},g\in\GL_2(\R)\}$$
be the conjugacy class of the morphism $h_0$.

An easy computation shows that the centralizer of $h_0$
in $\GL_2(\R)$ is the subgroup $D(\R)$.
We thus obtain an identification $X\cong\GL_2(\R)/D(\R)$.
There is also a left $\GL_2(\Z)$-action on $X$, and this gives us an 
interpretation of the space of geodesics as the quotient
$Y\cong \GL_2(\Z)\backslash X$.

The reader will probably ask now: what did we win in this translation?
The author's answer is: a strong analogy with the modular curve.

This allows us to view $Y$ as the moduli space of
triples $(M,F,\widetilde{F})$, where $M$ is a free $\Z$-module of rank $2$, and
$F,\widetilde{F}\subset M_\R$ are two lines in the corresponding real vector
space $M_\R:=M\otimes_\Z\R$.
To each $h\in X$, we associate a triple $(\Z^2,F_x,F_y)$ with $F_x$ the line
of weight $x$ for $h$ (i.e., $h(x,y)\cdot v=x\cdot v$ for all $v\in F_x$),
and $F_y$ the line of weight $y$.
It is helpful to think of the direct sum decomposition
$$M_\R=F\oplus \widetilde{F}$$
as an analog of the Hodge decomposition of the first
complex singular cohomology group of an elliptic curve $E/\C$:
$$H^1(E,\C)=H^1(E,\Oc_E)\oplus H^0(E,\Omega^1_E).$$

\subsection{Examples of bad Mumford-Tate groups}
In view of the analogy in the previous paragraph, we can ask:
for geodesics, what is the analog of the Mumford-Tate
group\footnote{A kind of analytic Galois group.} of elliptic curves?
Recall that an element $h\in X$ is a morphism of algebraic groups
$h:\G_{m,\R}^2\to \GL_{2,\R}$.

We suggest two possible analogs. The first is obtained by copying the
usual definition.

\begin{definition}
Let $h\in X$. We define the \emph{bad Mumford-Tate group} of $h$ to be
the smallest $\Q$-algebraic subgroup $\BMT(h)\subset \GL_{2,\Q}$ such that
$h((\R^*)^2)\subset \BMT(h)(\R)$.
\end{definition}

Let us test this definition on some examples.
It is clear that $\BMT(h_0)=D\cong \G_{m,\Q}^2$ is the group of rational
diagonal matrices, the maximal torus of $\GL_{2,\Q}$.

Let $u\in\R$ and
$g_u:=\left(\begin{smallmatrix}
1 & u\\
0 & 1
\end{smallmatrix}
\right)$ be the corresponding unipotent matrix.
Denote $h_u:=g_uh_0g_u^{-1}\in X$.
The morphism $h_u$ is given by the matrix
$\left(\begin{smallmatrix}
x & u(y-x)\\
0 & y
\end{smallmatrix}
\right).$
If we suppose that $u$ is not rational, then $\BMT(h_u)$ is the subgroup
$B\subset\GL_{2,\Q}$
of upper triangular matrices. Otherwise,
$\BMT(h_u)=g_uDg_u^{-1}\cong \G_{m,\Q}^2$ is a maximal torus.

Let $h'$ be the conjugate of $h_0$ by the matrix
$g':=\left(\begin{smallmatrix}
1 & 1\\
u & -u
\end{smallmatrix}
\right)$, with $u=e\in\R$. Then we get $h'\in X$ 
such that $\BMT(h')=\GL_{2,\Q}$.

Now let us consider a square free positive integer $d>1$.
Let
$g:=\left(\begin{smallmatrix}
1 & 1\\
\sqrt{d} & -\sqrt{d}
\end{smallmatrix}
\right)$.
If we conjugate $h_0$ by $g$, we obtain the matrix
$h'=\left(\begin{smallmatrix}
a & b\\
db & a
\end{smallmatrix}
\right)$ with $a=\frac{x+y}{2}$ and $b=\frac{x-y}{2\sqrt{d}}$.
The group of matrices of the form
$\left(\begin{smallmatrix}
a & b\\
db & a
\end{smallmatrix}
\right)$
with $a$ and $b$ rational is an algebraic torus over $\Q$,
conjugated over $\R$ to the maximal torus $D(\R)$.
It is clearly the Mumford-Tate group of $h'$.

\subsection{Good Mumford-Tate groups}
We now want to modify the definition of bad Mumford-Tate groups to obtain
``good'' ones.
The examples of last paragraph give us a quite precise idea of the different
kinds of bad Mumford-Tate groups that can appear. However, we would like to
have \emph{reductive} Mumford-Tate groups in order to have a closer
analogy with the case of elliptic curves.

\begin{definition}
\label{defin-MT-geodesiques}
Let $h\in X$. A \emph{Mumford-Tate group} for $h$ is a minimal reductive
subgroup of $\GL_{2,\Q}$ that contains $\BMT(h)$.
\end{definition}
With this new definition, a Mumford-Tate group is always reductive
whereas $\BMT$ may be, for example, the group $B$ of upper triangular
matrices.

\begin{lemma}
Let $h\in X$. Then there exists a unique Mumford-Tate group for $h$,
which we denote by $\MT(h)$.
\end{lemma}
\begin{proof}
The group $\BMT_\R$ contains a maximal torus of $\GL_{2,\R}$
(the image of $h$).
We thus have three possibilities, depending on the dimension
($2$,$3$ or $4$) of $\BMT$: $\BMT$ is a maximal torus,
a Borel subgroup or the whole of $\GL_{2,\Q}$.
If $\BMT$ is a maximal torus or $\GL_{2,\Q}$ then $\MT$ is clearly well
defined and equal to $\BMT$.
If $\BMT$ is a Borel subgroup, then the smallest reductive group that
contains it is the group $\GL_{2,\Q}$.
\end{proof}

We now arrive at the structure theorem for Mumford-Tate groups of geodesics.

\begin{theorem}
\label{classif-MT-geodesiques}
Let $h\in X$. The Mumford-Tate group of $h$ is of one of the following types:
\begin{enumerate}
\item $\MT(h)=\GL_2$, the corresponding geodesic is called MT-generic,
\item $\MT(h)=\Res_{E/\Q}\G_m$ with $E/\Q$ a real quadratic field,
we say that the corresponding geodesic is \emph{special} or has
\emph{real multiplication by $E$},
\item $\MT(h)=\G_{m,\Q}^2$, the geodesic is rational.
\end{enumerate}
\end{theorem}

This theorem follows from the the fact that all maximal tori in
$\GL_{2,\Q}$ are of the form $\Res_{E/\Q}\G_m$ with $E/\Q$
\'etale of dimension $2$. Such an algebra $E$ is either a real quadratic field,
isomorphic to $\Q^2$, or an imaginary quadratic field.
The imaginary quadratic
case can not appear because the two lines $F_x$ and $F_y$
associated to $h$ (see Subsection \ref{analogs-hodge-structures}) are real.

The properties of Mumford-Tate groups have a simple interpretation in terms
of dynamical properties of geodesics on the modular curve, as explained
to the author by Etienne Ghys.
The first case corresponds to non-closed geodesics,
the second to closed geodesics not homotopic to the cusps,
and the last corresponds to closed geodesics homotopic to the cusps.

\subsection{Special geodesics and class field theory}
\label{special-geodesics-CFT}
Let us look more closely at the second case
in Theorem \ref{classif-MT-geodesiques}, following Gauss' ideas.
So let $h\in X$ correspond to a special geodesic (i.e., one with real
multiplication by a real quadratic field $E$) and
let $T$ be its Mumford-Tate group. By definition,
$h$ is a morphism into $T_\R$.  The $T(\R)$-conjugacy class
of $h$ is reduced to $\{h\}$. The natural map
$(T,\{h\})\to (\GL_{2,\Q},X)$ induces a map between the double
coset spaces
$$
\xymatrix{
T(\Q)\backslash (\{h\}\times T(\A_f)/T(\widehat{\Z}))\ar[r]\ar@{=}[d] &
\GL_2(\Q)\backslash (X\times \GL_2(\A_f)/\GL_2(\widehat{\Z}))\ar@{=}[d]\\
\Pic(\Oc_E)\ar[r] & \GL_2(\Z)\backslash X=Y}
$$
whose image gives the space of geodesics with real multiplication by $E$.
The left term of this morphism is the ideal class group $\Pic(\Oc_E)$
of this real multiplication field.
The bottom arrow interprets $\Pic(\Oc_E)$ in terms of geodesics, a well
known result of Gauss.
The top arrow is an adelic formulation of this result which may admit a
generalization to higher rank spaces.

We want to arrive at a similar geodesic interpretation for the connected
component of the id\`ele class group $\pi_0(T(\Q)\backslash T(\A))$.
Recall that if $(G,H)$ is a pair consisting of
a reductive group over $\Q$ and a $G(\R)$-left space $H$, we denote by
$$
\Sh(G,H)=\limproj_K G(\Q)\backslash(H\times G(\A_f)/K)
$$
where $K\subset G(\A_f)$ runs over all compact open subgroups.
In \cite{De4}, 2.2.3, it is shown that
$$\pi_0(T(\Q)\backslash T(\A))=\Sh(T,\pi_0(T(\R))).$$
At this point, the morphism $(T,\pi_0(T(\R)))\to (\GL_2,X)$ induces a map
\begin{equation}
\label{injection-Shimura}
\Sh(T,\pi_0(T(\R)))\to \Sh(\GL_2,X).
\end{equation}
Unfortunately, the archimedian component $\pi_0(T(\R)):=T(\R)/T(\R)^+$
is killed under this mapping because $X\cong \GL_2(\R)/T(\R)$.
Therefore, one replaces $X$ by $X^\pm\cong \GL_2(\R)/T(\R)^+$
in morphism (\ref{injection-Shimura}) where $X^\pm$ is
the space of morphisms $h\in X$ \emph{oriented}
by the choice of orientations $s_x$, $s_y$ on the corresponding real lines
$F_x$, $F_y$.
This space projects naturally onto $X$.
The resulting map
$$\pi_0(T(\Q)\backslash T(\A))=\Sh(T,\pi_0(T(\R)))\to \Sh(\GL_2,X^\pm)$$
yields the geodesic interpretation of the connected component of the
id\`ele class group.

We are however far from a theory of real multiplication because
a ``natural'' rational structure on the space of geodesics,
analogous to the coordinates on the modular curve given
by the $j$ and $\Pc$ functions, is missing.

\section{Algebraic groups and moduli spaces of non-commutative tori}
\subsection{The moduli space of elliptic curves}
\label{famille-univ-classique}
This section recalls the classical construction of
the universal family of elliptic curves in terms of mixed Shimura varieties
\cite{Pink}, 10.7.
For classical Shimura varieties, we will refer to
\cite{De4}. The basic definition of mixed Shimura varieties can
be found in \cite{Pink}, 2.1, 3.1, or \cite{Mi3}[VI]. We only need in
this paper the notion of pre-Shimura datum defined in the Notations
Section.

Let $\Sb:=\Res_{\C/\R}\G_m$ be the Deligne torus. This is an $\R$-algebraic
group such that $\Sb(\R)=\C^*$.
Let $w:\G_{m,\R}\to \Sb$ be the weight morphism given by the natural inclusion
$\R^*\subset \C^*$.
Let
$\mu:\G_{m,\C}\to \Sb_\C\cong \G_{m,\C}\times \G_{m,\C}$ be the
\emph{Hodge morphism}
that sends $z$ to the pair $(z,1)$ and let $\widetilde{\mu}$ the morphism
that sends $z$ to the pair $(1,z)$. We call $\widetilde{\mu}$ the
\emph{anti-Hodge morphism}.

Let $V$ be an $\R$-vector space.
Recall that a representation $h:\Sb_\C\to \GL(V_\C)$ in the $\C$-vector space
$V_\C:=V\otimes_\R\C$
gives an ascending filtration, the so called \emph{weight filtration},
$W_\bullet V_\C$ given by the
cocharacter $h\circ w_\C$ and 
a descending filtration,
the so called \emph{Hodge filtration}, $F^\bullet V_\C$ given by the
cocharacter $h\circ \mu$. We will also be interested by another descending,
that we call the \emph{anti-Hodge filtration},
$\widetilde{F}^\bullet V_\C$ given
by the cocharacter $h\circ \widetilde{\mu}$. Note that if $h$
is not defined over $\R$, then the anti-Hodge filtration is
not the complex conjugate filtration of $F^\bullet V_\C$.

Let $\Hb^{\pm}:=\{\tau\in\C,\tau\notin\R\}$ be the Poincar\'e double half
plane. This space can be identified with the $\GL_2(\R)$-conjugacy class of the 
morphism of real algebraic groups
$h_i:\Sb\to \GL_{2,\R}$
that maps $z=a+ib\in\C^*=\Sb(\R)$ to the matrix
$\left(\begin{smallmatrix}
a & b\\
-b & a
\end{smallmatrix}\right).$
The map between this conjugacy class and $\Hb^{\pm}\subset \Pb^1(\C)$ is
given by associating to $h:\Sb\to\GL_{2,\R}$ the
line $F(h):=\Ker(h_\C\circ\mu(z)-z\cdot\id)\subset \C^2$.

\begin{definition}
The pair $(\GL_2,\Hb^{\pm})$ is called the \emph{classical modular
Shimura datum}.
\end{definition}

Let $K$ be the compact open subgroup $\GL_2(\hat{\Z})$ of $\GL_2(\A_f)$.
The associated Shimura variety of level $K$ is by definition
$$
\begin{array}{rl}
\Mcal=\Sh_K(\GL_{2,\Q},\Hb^{\pm})&
=\GL_2(\Q)\backslash (\Hb^{\pm}\times \GL_2(\A_f)/K)\\
& =\GL_2(\Z)\backslash \Hb^{\pm}\\
& =\PGL_2(\Z)\backslash\Hb,
\end{array}
$$
that is to say the classical modular curve.

Now we recall the construction of the universal family $\Ec\to \Mcal$ of
elliptic curves.

Let $P$ be the group scheme $V\rtimes\GL_2$ with $V:=\G_a^2$ the standard 
representation\footnote{i.e $V(\Z)=\Z^2$} of $\GL_2$. Fix a rational splitting
$s:\GL_2\to P$ of the natural projection map $\pi:P\to \GL_2$,
for example the map $g\mapsto (0,g)\in V\rtimes\GL_2$.
All such splittings are conjugates under $V(\Q)\subset P(\Q)$.
Define
$h':\Sb_\C\to P_\C$ by $h':=s_\C\circ h_{i,\C}$.
We will denote by $\Hb'$ the $P(\R)$-orbit of $h'$ in $\Hom(\Sb_\C,P_\C)$
for the conjugation action.

\begin{definition}
The pair $(P,\Hb')$ is called the \emph{pre-Shimura datum of the universal
elliptic curve}.
\end{definition}

\begin{lemma}
\label{indep-section}
The space $\Hb'$ does not depend on the choice of the splitting $s$.
\end{lemma}
\begin{proof}
This follows from the fact that all such splittings are conjugates under
$V(\Q)$ and that $\Hb'$ is a $P(\R)$-orbit.
\end{proof}

Let $K^P:=P(\hat{\Z})\subset P(\A_f)$.
The mixed Shimura variety associated to the data $(P,\Hb')$ and $K^P$ is, by
definition,
$$\Ec=\Sh_{K^P}(P,\Hb')=
P(\Q)\backslash (\Hb'\times P(\A_f)/K^P)=P(\Z)\backslash \Hb'.$$

The natural projection map $\pi:P\to \GL_2$ induces a projection morphism
of Shimura data
$(P,\Hb')\to (\GL_2,\Hb^{\pm})$ and
a projection map
$$\Ec\to \Mcal$$
between the corresponding mixed Shimura varieties.

Theorem 10.10 of \cite{Pink} shows that this map (up to the choice of
finer level structures $K\subset \GL_2(\A_f)$ and $K^P\subset P(\A_f)$) gives
the universal family of elliptic curves.

\begin{remark}
We want to stress here that the fiber of
$\Hb'\to \Hb^{\pm}$ over some $h:\Sb\to \GL_{2,\R}$ is given by
the one dimensional $\C$-vector space $\C^2/F(h)\cong V(\R)$.
This will be useful in the next section.
\end{remark}

\subsection{The universal family of geodesics}
\label{famille-univ-geodesiques}
We will now use the same ideas to construct, for the space of
geodesics, an analog of the universal family of elliptic curves.

Let $h^Q_0:\G_{m,\R}^2\to \GL_{2,\R}^2$ be the morphism that sends the pair 
$(x,y)\in(\R^*)^2$ to the pair of matrices
$\left
(\left(\begin{smallmatrix}
x & 0\\
0 & 1
\end{smallmatrix}\right),
(\left(\begin{smallmatrix}
1 & 0\\
0 & y
\end{smallmatrix}\right)\right).$
Let $Z_Q:=\GL_2(\R)\cdot h_0^Q$ be the $\GL_2(\R)$-conjugacy class
of $h_0^Q$ in $\Hom(\G_{m,\R}^2,\GL_{2,\R}^2)$. The multiplication
map $m:\GL_2^2\to \GL_2$ (which is not a group homomorphism) induces a natural
$\GL_2(\R)$-equivariant bijection
$$\begin{array}{ccc}
Z_Q & \overset{m}{\to} & X\\
\left(g\left(\begin{smallmatrix}
x & 0\\
0 & 1
\end{smallmatrix}\right)g^{-1},
g\left(\begin{smallmatrix}
1 & 0\\
0 & y
\end{smallmatrix}\right)g^{-1}\right)
& \mapsto & 
g\left(\begin{smallmatrix}
x & 0\\
0 & y
\end{smallmatrix}\right)g^{-1}.
\end{array}
$$
In other words, we have constructed an isomorphism
$m:(\GL_2,Z_Q)\to (\GL_2,X)$ of pre-Shimura data.

Let $V:=\G_a^2$ be the standard representation\footnote{i.e., $V(\Z)=\Z^2$.}
of $\GL_2$. Let $Q'$ be the group scheme $V^2\rtimes \GL_2^2$ and
$Q$ be the group scheme $V^2\rtimes \GL_2$.
We will also denote by $h^Q_0:\G_{m,\R}^2\to Q'_\R$ the morphism obtained
by composing $h^Q_0$ with a rational section of the natural projection
$Q'\to \GL_2^2$. Such a section is unique up to an element of $V^2(\Q)$.
Let $Y_Q:=Q(\R)\cdot h^Q_0\subset \Hom(\G_{m,\R}^2,Q'_\R)$ be the
$Q(\R)$-conjugacy class of $h^Q_0$. It does not depend on the chosen
section because $V^2\subset Q$.

\begin{definition}
We call the pair $(Q,Y_Q)$ the
\emph{pre-Shimura datum of the universal family of geodesics}.
\end{definition}
The next lemma will explain this definition.

Let $\pi':Q'\to\GL_2^2$ be the natural projection. This projection induces
a natural map $Y_Q\to Z_Q$ that is compatible with the projection
$\pi:Q\to\GL_2$. This yields a morphism of pre-Shimura data
$\pi:(Q,Y_Q)\to (\GL_2,Z_Q)$.
If we compose this morphism with $m$, we get
a natural morphism of pre-Shimura data
$$(Q,Y_Q)\to (\GL_2,X)$$
which is in fact the quotient map by the additive group $V^2$.

The Shimura fibered space $\Sh(Q,Y_Q)\to \Sh(\GL_2,X)$ can be considered
as a universal family of geodesics because of the following lemma.

\begin{lemma}
\label{lemme-famille-univ-geodesiques}
The fiber of the projection
$\Sh_{Q(\widehat{\Z})}(Q,Y_Q)\to \Sh_{\GL_2(\widehat{\Z})}(\GL_2,X)$
over a point $[(V,F_x,F_y)]$ of the space of geodesics is the space
$\Z^2\backslash \R^2/F_x\times \Z^2\backslash \R^2/F_y$, product
of the two leaves spaces of the corresponding linear foliations on
the torus $\Z^2\backslash \R^2=V(\Z)\backslash V(\R)$.
\end{lemma}
\begin{proof}
We can embed $Q'$ in $\GL_3^2$ by considering pairs of matrices of
the form
$$\left(\left(\begin{smallmatrix}
A_1 & v_1\\
0 & 1
\end{smallmatrix}\right),
\left(\begin{smallmatrix}
A_2 & v_2\\
0 & 1
\end{smallmatrix}\right)\right)$$
with $A_1,A_2\in \GL_2$ and $v_1,v_2\in V$. The morphism
$h_0^Q$ is now given by the matrix
$\left(\left(\begin{smallmatrix}
Z_x & 0\\
0 & 1
\end{smallmatrix}\right),
\left(\begin{smallmatrix}
Z_y & 0\\
0 & 1
\end{smallmatrix}\right)\right)$
with
$Z_x=\left(\begin{smallmatrix}
x & 0\\
0 & 1
\end{smallmatrix}\right)$,
$Z_y=\left(\begin{smallmatrix}
1 & 0\\
0 & y
\end{smallmatrix}\right)$. The $Q(\R)$-conjugacy class $Y_Q$ of
$h_0^Q$ is given by pairs of matrices of the form
$$
(M_x,M_y):=\left(\left(\begin{smallmatrix}
AZ_xA^{-1} & (I-AZ_xA^{-1})v_1\\
0 & 1
\end{smallmatrix}\right),
\left(\begin{smallmatrix}
AZ_yA^{-1} & (I-AZ_yA^{-1})v_2\\
0 & 1
\end{smallmatrix}\right)\right)
$$
with $A\in\GL_2(\R)$. The projection
$Y_Q\to Z_Q$ sends such a matrix to the pair
$(AZ_xA^{-1},AZ_yA^{-1})$. The action of $(v_1',v_2')\in V^2(\R)$ by conjugation
on a pair $(M_x,M_y)$ as above gives
$$\left(\left(\begin{smallmatrix}
AZ_xA^{-1} & (I-AZ_xA^{-1})(v_1+v_1')\\
0 & 1
\end{smallmatrix}\right),
\left(\begin{smallmatrix}
AZ_yA^{-1} & (I-AZ_yA^{-1})(v_2+v_2')\\
0 & 1
\end{smallmatrix}\right)\right).
$$
In this way, we see directly that the fibre $Y_{Q,h}$ of this projection
over some $h\in Z_Q\cong X$ is a $V^2(\R)$-homogeneous space.
The stabilizer of any
point of $Y_{Q,h}$ for this $V^2(\R)$-action is the sum
$F_y\oplus F_x\subset V(\R)\oplus V(\R)$. This shows that
$Y_{Q,h}\cong V(\R)/F_y\times V(\R)/F_x$. Since the projection is equivariant
with respect to the projection $Q\to \GL_2$, the fibre of the projection
$\Sh_{Q(\widehat{\Z})}(Q,Y_Q)\to \Sh_{\GL_2(\widehat{\Z})}(\GL_2,X)$
at the point $[\GL_2(\Z)\cdot h]$
is given by the space
$\Z^2\backslash \R^2/F_x\times \Z^2\backslash \R^2/F_y$.
\end{proof}

\subsection{The moduli problem for non-commutative tori}
In this section, we want to understand Manin's point of view of the
moduli space of non-commutative tori (in terms of pseudo-lattices)
in the spirit of Hodge structures. To this end, we introduce the notion
of pre-lilac.

\begin{definition}
A (rank $2$) \emph{pre-lilac}\footnote{Abbreviation for \textbf{LI}ne
in a \textbf{LA}tti\textbf{C}e.}
is a pair $(M,F)$ of a free $\Z$-module of rank two and
a real line $F\subset M_\R$.
A \emph{morphism of pre-lilacs} $(M_1,F_1)\to (M_2,F_2)$ is
a morphism of Abelian groups $f:M_1\to M_2$ such
that $f_\R(F_1)\subset F_2$.
\end{definition}

This notion is equivalent to Manin's notion of pseudo-lattice
with weak morphisms (see \cite{Manin3}) but it is easier for us to formulate
our results in terms of lilacs because of their analogy with complex
structures.

From a pre-lilac $(M,F)$, one can construct a non-commutative
algebra $\Ac(M,F)$, called the \emph{Kroneker foliation algebra}
$C^\infty(M\backslash M_\R)\rtimes F$.
We will call such an algebra a \emph{non-commutative torus} because the
choice of an element $e$ of a basis for $M$ allows one to construct a Morita
equivalent algebra
$\Tc(M,F,e)=C^\infty(\Z.e\backslash M_\R/F)\rtimes [\Z.e\backslash M]$
which is an irrational rotation algebra, i.e., a non-commutative torus
in the usual sense. The Morita equivalence
$\Ac(M,F)\sim\Tc(M,F,e)$ follows from
\cite{ENCG}, corollary 12.20\footnote{See also \cite{ENCG},
theorem 12.17 for a choice-of-basis free proof of this fact.}.
We are interested in the
set of isomorphism classes of Kronecker foliation algebras.
Let $\nct$ be the category of such algebras with
$*$-isomorphisms as morphisms. Let $\prelilacs$ be the category
of pre-lilacs with isomorphisms as morphisms.
The assignment $(M,F)\mapsto \Ac(M,F)$ gives a functor
$$
\begin{array}{cccc}
T: & \prelilacs & \to   & \nct\\

   & (M,F) &\mapsto& \Ac(M,F).
\end{array}
$$

For an object $A$ of $\nct$, denote by
$S:\HC_2(A)\to\HC_0(A)$ the periodicity map in cyclic homology,
by $\ch:K_0(A)\to HC_2(A)$ the Chern character
defined in \cite{Loday}, and by $\ch_\R:K_0(A)\otimes_\Z\R\to HC_2(A)$
the corresponding map of $\R$-vector spaces.

There is also a functor in the other direction
$$
\begin{array}{cccc}
L: & \nct 	& \to   & \prelilacs\\
   & A   	&\mapsto& (K_0(A),
                           \ch_\R^{-1}(\Ker(S))\subset K_0(A)\otimes_\Z\R)
\end{array}
$$
The facts that $\ch$ induces an isomorphism $K_0(A)\otimes_\Z\C\to HC_2(A)$,
and that the filtration by $\Ker(S)$ on cyclic homology is real with
respect to the real structure given by $K$-theory can be deduced from
the explicit calculations of Lemma 54 of \cite{Connes3}.
Notice that this functor was already present in the paper \cite{Connes1}.

\begin{lemma}
Two Kronecker foliation algebras $\Ac_1$ and $\Ac_2$ are Morita equivalent
if and only if they are isomorphic as $*$-algebras.
\end{lemma}
\begin{proof}
This comes from the fact that for every objects $\Ac(M,F)\in \nct$, one has
$L(\Ac(M,F))\cong (M,F)$.
To prove this, we can look at the Morita equivalent algebra $\Tc(M,F,e)$, for
$e$ an element of a basis for $M$ and use the explicit calculation of its
Chern character in \cite{Connes3}.
We thus obtain that for all object $\Ac\in\nct$, $T\circ L(\Ac)\cong \Ac$.
To finish the proof, if the two Kronecker algebras are Morita equivalent,
then Morita invariance of cyclic cohomology, K-theory and Chern character
implies $L(\Ac_1)\cong L(\Ac_2)$ and this implies $\Ac_1\cong\Ac_2$.
\end{proof}

The relation with Manin's definition of pseudo-lattices associated to
non-commutative tori is the following.
The long exact sequence of cyclic homology gives
$$
HH_2(A)\overset{I}{\to} HC_2(A)\overset{S}{\to}
HC_0(A)\overset{B}{\to} HH_1(A)
$$
and we know that $S$ is surjective in this case
(by explicit calculation given in
\cite{Connes3}), so it gives an isomorphism
$HC_2(A)/\Ker(S)\overset{S}{\to} HC_0(A)$.
Since the Chern character is compatible with $S$, the natural map
$S\circ\ch: K_0(A)\to HC_0(A)$ induced by $\ch:K_0(A)\to HC_2(A)$ is
equal to the Chern character
$\ch:K_0(A)\to HC_0(A)$.
The pseudo-lattices that Manin considers are given by this Chern character.
So the functor that associates to a pre-lilac $(M,F)$ the pseudo-lattice
$(M,M_\C/F_\C)$ has a natural interpretation in cyclic homology as\
the association
$(K_0(A),\ch_\R^{-1}(\Ker(S)))\leadsto (K_0(A),HC_0(A)=HC_2(A)/\Ker(S))$.

The two functors $T$ and $L$ naturally identify the \emph{set} of isomorphism
classes of pre-lilacs and the \emph{set} of isomorphism classes of
non-commutative tori.

Fix a free $\Z$-module $M$ of rank two. The projective space $\Pb(M_\R)$ over
$M_\R$ gives a parameter space for lines $F'\subset M_\R$. Two such lines
correspond to isomorphic pre-lilacs if and only if they are exchanged by
some $g\in \GL(M)$. So the set of isomorphism classes of pre-lilacs is
$\GL(M)\backslash P(M_\R)\cong \GL_2(\Z)\backslash \Pb^1(\R)$. As said before,
this is also the moduli set for non-commutative tori.

We also remark that there is a natural projection
$$
\begin{array}{ccc}
\GL_2(\Z)\backslash X 	&\to    & \GL_2(\Z)\backslash\Pb^1(\R)\\
(M,F_x,F_y)		&\mapsto& (M,F_x)
\end{array}
$$
from the space of geodesics to the moduli space of pre-lilacs.
In non-commutative geometry, the left hand side of this projection
can be interpreted as the moduli space for triples
$$(\Ac_1,\Ac_2,\psi:K_0(\Ac_1)\overset{\sim}{\to} K_0(\Ac_2))$$
consisting of two Kronecker foliation algebras
and an isomorphism $\psi$ between their $K_0$ groups, such that if $F_1$ and
$F_2$ are the real lines in $K$-theory constructed in this section,
$\psi_\R(F_1)\oplus F_2=K_0(\Ac_2)_\R$.

\begin{remark}
As we saw in Section \ref{special-geodesics-CFT},
the study of nontrivial level structures on geodesics (and thus on
non-commutative tori) is
also interesting because of class field theory.
We will denote by $\Irrat(X)^\pm\subset X^\pm$ the space of pairs
of irrational oriented lines in $\R^2$.
Let $N>1$ be an integer, and let $K_N$ be the group
defined by the exact sequence
$1\to K_N\to \GL_2(\widehat{\Z})\to \GL_2(\Z/N\Z)\to 1.$
Then the Shimura space
$\Sh_{K_N}(\GL_2,\Irrat(X)^\pm)$
is the moduli space of tuples
$$(M,F_x,F_y,s_x,s_y,\phi:M\otimes_\Z\Z/N\Z\overset{\sim}{\to} (\Z/N\Z)^2)$$
consisting of a free $\Z$-module $M$ of rank $2$, two irrational lines
$F_x$, $F_y$ in the underlying real vector space, equiped with two
orientations $s_x$, $s_y$, and a level structure $\phi$.
Perhaps this space has a moduli interpretation in terms of tuples
$$
(\Ac_1,\Ac_2,
\psi:K_0(\Ac_1)\overset{\sim}{\to} K_0(\Ac_2),
\phi:K_0(\Ac_1)\otimes_\Z \Z/N\Z\overset{\sim}{\to} (\Z/N\Z)^2).
$$
The orientation on the real lines could be given by the image of the positive
cone in $K$-theory. The author's knowledge is, however, not sharp enough
to be sure of this non-commutative moduli interpretation.
\end{remark}

\begin{remark}
Polishchuk's remarkable work \cite{Poli2} on the relation of analytic
non-commutative tori with nonstandard t-structures on derived categories
of coherent sheaves on usual elliptic curves seems to be
promising, because it allows one to give some rationality
properties to the objects we have on hand.
For example, if $E$ is an elliptic curve over
$\Q$ and we fix some $t$-structure on its associated analytic curve
given by a line with real quadratic slope on $K_0(E)\otimes_\Z\R$,
we can ask questions about the rationality of this $t$-structure.
The Algebraic Proj construction in \cite{Poli1}
for real multiplication non-commutative tori
also allows one to ask rationality questions.
\end{remark}

\subsection{The universal non-commutative torus}
For at least four reasons (esthetic symmetry,
class field theory, Mumford-Tate groups
and the complications that appear in this paragraph),
it seems to be more natural to study the moduli space for \emph{pairs}
of non-commutative tori as above (which is also the space of geodesics)
rather than the moduli space of solitary non-commutative tori.
However, we still want to construct a universal non-commutative
torus because it permits to understand the relation with Tate mixed Hodge
structures.

Let $h_0'$ be the morphism of algebraic groups given by
$$
\begin{array}{cccc}
h_0': & \G_{m,\R}^2 & \to & \GL_{2,\R}\\
& (x,y) & \mapsto &
\left(\begin{smallmatrix}
xy & 0\\
0 & 1
\end{smallmatrix}\right)
\end{array}.
$$

\begin{aside}[rational boundary component]
If we denote by $U$ the group of unipotent upper triangular
matrices and by $P_1$ the group of matrices of the form
$\left(\begin{smallmatrix}
* & *\\
0 & 1
\end{smallmatrix}\right)$ in $\GL_{2,\Q}$,
then the $U(\C)$-conjugacy class $Y_1'$ of $h_{0,\C}'$ is
identified with $U(\C)=\C$. Let $Y_1$ be $Y_1'\times\{\pm 1\}$.
Pink calls the pair $(P_1,Y_1)$ a \emph{rational boundary component} of
$(\GL_2,\Hb^\pm)$.
These kind of (mixed) Shimura data appear in the
toroidal compactification of Shimura varieties. They are parameter
spaces for mixed Hodge structures. In our particular case, the mixed
Hodge structures are extensions
$$0\to \Z(1)\to M\to \Z\to 0$$
that can be described algebraically by $1$-motives of the form
$[\Z\to \G_{m,\C}]$ (see \cite{De11} section 10 for a definition of 
$1$-motives). Level structures on these one-motives are strongly
related to roots of unity in $\C$, i.e., to generators of $\Q^{ab}$.
In some sense (see \cite{Pink}, 10.15 to give a precise
meaning to this affirmation),
we can say that $\Sh(P_1,Y_1)$ is a universal family
over the moduli space $\Sh(\G_{m,\Q},\{\pm 1\})$ of primitive roots of unity.
\end{aside}

We will denote by $X'$ the $\GL_2(\R)$-conjugacy class of $h_0'$.
Recall that $X'$ is the set $\{(F_{xy},F_0)\}$ of pairs of distinct lines
in $\R^2$ of respective weights $0$ and $xy$. There is
a natural projection map $X'\to \Pb^1(\R)$
given by $(F_{xy},F_0)\mapsto F_0$.

Recall from Section \ref{famille-univ-classique} that
$P$ is the group scheme $V\rtimes\GL_2$, with $V:=\G_a^2$
the standard representation of $\GL_2$. As usual, we also
denote by $h_0':\G_{m,\R}^2\to P_\R$ the morphism obtained
by composing $h_0'$ with a rational section of the natural projection
$P\to \GL_2$. Such a section is unique up to an element of $V(\Q)$.
Let $Y_P$ be the $P(\R)$-conjugacy class of $h_0'$. It is independent
of the choice of the section of $P\to \GL_2$ because $V\subset P$.

Let $V_P=\Q^3$ be the standard representation of $P$ given by an embedding
$P\hookrightarrow \GL_3$.
For $h\in Y_P$,
we let $F^0(h):=\{v\in V_{P,\R}|h(x,y)\cdot v=v\}$.
On $Y_P$, we have the equivalence relation:
$$h\sim h' \Leftrightarrow F^0(h)V'=F^0(h')V'.$$
Let $\overline{Y_P}=Y_P/\sim$ be the corresponding quotient.

\begin{definition}
We call the pair $(P,\overline{Y_P})$ the
\emph{pre-Shimura datum of the universal family of non-commutative tori}.
\end{definition}
The following lemma shows us that the Shimura fibered space
$\Sh(P,\overline{Y_P})\to \Sh(\GL_2,\Pb^1(\R))$
can be considered as a \emph{universal family of non-commutative tori}.
\begin{lemma}
The fiber of the projection
$$
\Sh_{P(\widehat{\Z})}(P,Y_P)\to
\Sh_{\GL_2(\widehat{\Z})}(\GL_2,\Pb^1(\R))=
\PGL_2(\Z)\backslash \Pb^1(\R)
$$
over a point $(M,F_0)$ of the space of pre-lilacs is the space
$M\backslash M_\R/F_0$ of leaves of the corresponding foliation on
the two-torus $M\backslash M_\R$.
\end{lemma}
\begin{proof}
The proof is essentially the same as the one of Lemma
\ref{lemme-famille-univ-geodesiques}. The additional fact to check
is that forgetting the $F_{xy}$ part of the filtration is
compatible with the quotient, which follows from the definition.
\end{proof}

\section{Two higher dimensional examples}
\subsection{Hilbert modular varieties}
Let $E$ be a totally real number field with ring of integer $\Oc_E$,
let $I:=\Hom(E,\R)$ and let $n:=\card(I)$.
Denote by
$G$ the group scheme $\Res_{\Oc_E/\Z} \GL_2$.
We then have
$G_\R=\prod_{\iota:E\to \R}\G_{\R,\iota}\cong \prod_{\iota:E\to\R} \GL_{2,\R}$.

Let $h:=\prod_{\iota:E\to\R} h_i:\Sb\to\G_\R$ with
$h_i:\Sb\to \GL_{2,\R}$
the map that sends $z=a+ib\in\C^*=\Sb(\R)$ to the matrix
$\left(\begin{smallmatrix}
a & b\\
-b & a
\end{smallmatrix}\right).$
Let $X$ be the $G(\R)$-conjugacy class of this morphism. We have an isomorphism
$X\cong \prod_{\iota:E\to \R} \Hb^{\pm}$.

\begin{definition}
The Shimura datum $(G,X)$ is called the \emph{Hilbert modular Shimura datum}.
\end{definition}
Let $h(E)$ be the Hilbert class number of $E$.
Let $K$ be the compact open subgroup $G(\hat{\Z})$ of $G(\A_f)$.
The associated Shimura variety is by definition
$$
\begin{array}{rl}
\Mcal=\Sh_K(G,X)&
=G(\Q)\backslash (X\times G(\A_f)/K)\\
& =\GL_2(\Oc_E)\backslash X\textrm{ if $h(E)=1$,}\\
&\cong\GL_2(\Oc_E)\backslash(\prod_{\iota:E\to \R}\Hb^\pm)\textrm{ if $h(E)=1$,}
\end{array}
$$
i.e., the Hilbert modular variety.
Now let $B'$ be the subgroup of upper triangular matrices in $\GL_{2,\Oc_E}$
and denote by $B$ the group scheme $\Res_{\Oc_E/\Z}B'$.

\begin{aside}[rational boundary component]
The group scheme $B$ is a maximal parabolic subgroup of $G$ and corresponds to
a rational boundary component $(P_1,X_1)$ of $(G,X)$ as in \cite{Pink}, 4.11.
The canonical model of the associated mixed Shimura variety is a moduli space
defined over $\Q$ for $1$-motives with additional structures.
\end{aside}

Let $h_0^H:\G_{m,\R}^2\to G_\R$ be the morphism given on each simple component
of $G_\R$ by $h_0:(x,y)\in(\R^*)^2\mapsto \diag(x,y):=
\left(\begin{smallmatrix}
x & 0\\
0 & y
\end{smallmatrix}\right)$. Let $X_H$ be the $G(\R)$-conjugacy class of $h_0^H$.

\begin{definition}
We will call the pair $(G,X_H)$ the
\emph{pre-Shimura datum of the moduli space of Hilbert lilacs}.
\end{definition}

The corresponding Shimura space $\Sh_{G(\hat{\Z})}(G,X_H)$
is a moduli space for tuples $(M,F_x,F_y,i)$ consisting of a rank $2n$ free
$\Z$-module $M$, equipped with a decomposition
$$M_\R=F_x\oplus F_y$$
of its underlying real vector space in two $n$-dimensional subspaces,
and with a morphism $i:E\hookrightarrow \End(M_\Q,F_x,F_y)$ (i.e., a morphism
compatible with the decomposition). These objects are called
\emph{Hilbert lilacs}.

\begin{example}
Let $E:=\Q(\sqrt{2})$ and $F:=\Q(\sqrt{2},\sqrt{3})$. For each embedding
$\iota$ in $\Hom(E,\R)$,
choose one embedding $c_\iota$ over it in $\Hom(F,\R)$. Such a choice
is called an \emph{RM type} for $F/E$.
Equip $M:=\Oc_F$
with the ``Hodge decomposition''
$$M_\R=F_x\oplus F_y$$
where $F_x:=\oplus_\iota F_{c_\iota}\cong \R^2$, and $F_y$ is the
other component
in the natural decomposition of $M_\R$. Then the lilac $(M,F_x,F_y)$ gives
a point in the space of Hilbert lilacs corresponding to $E$, and this point
has Mumford-Tate group contained in
$\Res_{F/\Q}\G_m\subset G=\Res_{E/\Q}\GL_2$.
It is called a \emph{special point} of this space.
\end{example}

Let $h_0^{HQ}$ be the product
morphism $\prod_{\iota:E\to\R} h_0^Q:\G_{m,\R}^2\to G_\R^2$, where
$h^Q_0:\G_{m,\R}^2\to \GL_{2,\R}^2$ is the morphism that sends the pair 
$(x,y)\in(\R^*)^2$ to the pair of matrices
$\left
(\left(\begin{smallmatrix}
x & 0\\
0 & 1
\end{smallmatrix}\right),
(\left(\begin{smallmatrix}
1 & 0\\
0 & y
\end{smallmatrix}\right)\right).$

Let \ $Z_{HQ}:=G(\R)\cdot h_0^{HQ}$ \ be the \ $G(\R)$-conjugacy class
of \ $h_0^{HQ}$ \ in $\Hom(\G_{m,\R}^2,G_\R^2)$. The multiplication
map $m:G^2\to G$ (which is not a group homomorphism) induces, 
as in Subsection \ref{famille-univ-geodesiques}, a natural
isomorphism
$m:(G,Z_{HQ})\to (G,X_H)$ of pre-Shimura data.

Let $V:=\Res_{\Oc_E/\Z}\G_{a,E}^2$
be the standard representation\footnote{i.e., $V(\Z)=\Oc_E^2$.}
of $G$. Let $Q'$ be the group scheme $V^2\rtimes G^2$ and
let $Q$ be the group scheme $V^2\rtimes G$.
We will also denote by $h^{HQ}_0:\G_{m,\R}^2\to Q'_\R$ the morphism obtained
by composition with a rational section of the natural projection
$Q'\to G^2$. Such a section is unique up to an element of $V^2(\Q)$.

Let $Y_{HQ}:=Q(\R)\cdot h^{HQ}_0\subset \Hom(\G_{m,\R}^2,Q'_\R)$ be the
$Q(\R)$-conjugacy class of $h^{HQ}_0$. It does not depend on the chosen
section because $V^2\subset Q$.

\begin{definition}
We call the pair $(Q,Y_{HQ})$ the
\emph{pre-Shimura datum of the universal family of Hilbert lilacs}.
\end{definition}

The next lemma will explain this definition.
Let $\pi':Q'\to G^2$ be the natural projection. This projection induces
a natural map $Y_{HQ}\to Z_{HQ}$ that is compatible with the natural projection
$\pi:Q\to G$. This yields a morphism of pre-Shimura data
$\pi:(Q,Y_{HQ})\to (G,Z_{HQ})$.
If we compose this morphism with $m$, we get
a natural morphism of pre-Shimura data
$$(Q,Y_{HQ})\to (G,X_H)$$
which is in fact the quotient map by the additive group $V^2$.

The Shimura fibered space $\Sh(Q,Y_{HQ})\to \Sh(G,X_H)$ can be considered
as a universal family of Hilbert lilacs because of the following lemma.

\begin{lemma}
The fiber of the projection
$\Sh_{Q(\widehat{\Z})}(Q,Y_Q)\to \Sh_{G(\widehat{\Z})}(G,X)$
over a point $[(V,F_x,F_y,i)]$ of the space of Hilbert lilacs is the space
$V\backslash V_\R/F_x\times V\backslash V_\R/F_y$, product
of the two leaves spaces of the corresponding linear foliations on
the torus $V\backslash V_\R$.
\end{lemma}
\begin{proof}
The proof is essentially the same as in Lemma 
\ref{lemme-famille-univ-geodesiques}.
\end{proof}

\begin{remark}
The centralizer $C_{G(\R)}(h_0^H)$ of our basis morphism
is isomorphic to the maximal torus $T(\R)$ of $G(\R)$.
We define the
Mumford-Tate group of some $h\in X_H$ as in Definition 
\ref{defin-MT-geodesiques}
as a reductive envelope,
defined by Andr\'e, Kahn and O'Sullivan, \cite{Andre-Kahn}.
Such a reductive envelope is well defined up to the centralizer of the
enveloped group. We will now change a little bit our description of
$X_H$ in order to have a well defined Mumford-Tate group for all points
in this (and other) spaces. These groups could also be of some
interest in the study of dynamical properties of the corresponding
foliations, as suggested to the author by Yves Andr\'e and Etienne Ghys.
\end{remark}

Let $D$ be the basic maximal torus of $G$, i.e.,
$D:=\Res_{E/\Q}\G_m^2$. Let $\Db:=D_\R$ be the corresponding
real algebraic group, and let $h_{0,\Rc}^H:\Db\to G_\R$
be the natural inclusion. Let $\Rc_H$ be the $G(\R)$-conjugacy class
of $h_{0,\Rc}^H$. The inclusion\footnote{Induced by the rational inclusion
$\G_{m,\Q}^2\subset D$.}
$\G_{m,\R}^2\subset \Db$ induces a
$G(\R)$-equivariant bijection
$\Rc_H\to X_H$, i.e., an isomorphism of pre-Shimura data
$$(G,\Rc_H)\to (G,X_H).$$

\begin{definition}
Let $h\in \Rc_H$. We define the \emph{bad Mumford-Tate group} of $h$ to be
the smallest $\Q$-algebraic subgroup $\BMT(h)\subset G$ such that
$\im(h)\subset \BMT(h)(\R)$.
A \emph{Mumford-Tate group} for $h$ is a minimal reductive
subgroup of $G$ that contains $\BMT(h)$.
\end{definition}

\begin{lemma}
Let $h\in \Rc_H$. There exists a unique Mumford-Tate group for $h$.
It will be denoted by $\MT(h)$.
\end{lemma}
\begin{proof}
The bad Mumford-Tate group of $h$ contains a maximal torus $T$ of $G$
(because the image of $h$ is a maximal torus over $\R$).
We know from \cite{SGA3III}, XIX, 2.8 (or other classical references on
algebraic groups) that the centralizer $C_G(T)$ of $T$ in $G$
is $T$ itself, because $G$ is a reductive group.
We have an inclusion of centralizers $C_G(\BMT(h))\subset C_G(T)=T$,
and we already know that $T\subset \BMT(h)$. So we have
$C_G(\BMT(h))\subset \BMT(h)$, and the fact that a reductive envelope
of a group is well defined up to the centralizer of the group implies
that, in this case, the Mumford-Tate group is well-defined.
\end{proof}

\begin{remark}
This result gives a motivation to study the pair
$(G,\Rc_H)$ itself. This pair is called the
\emph{shore datum of the Hilbert modular datum $(G,X)$.}
Another good reason to use a bigger basis group $\Db$ for morphisms
is given by the need, in Manin's real multiplication program,
to take into account archimedian places in
class field theory of totally real fields.
This has been investigated in Section \ref{special-geodesics-CFT} in
the case $E=\Q$ and in
\cite{expose-rivages} for other quite general examples.
\end{remark}

\subsection{The moduli space of Abelian surfaces}
We first recall the construction of the moduli space of Abelian surfaces and
then construct one of its irrational boundaries. As we will see, this gives
a non-totally degenerate example, that is in some sense more general than
the case of Hilbert moduli spaces.

Let $V=\Z^4$ and $\psi:V\times V\to\Z$ be the standard symplectic form 
given by the matrix
$
J:=\left(\begin{smallmatrix}
0   & I_2\\
-I_2 & 0
\end{smallmatrix}\right)
$
with $I_2\in\mathrm{M}_2(\Z)$ the identity matrix.
Let $G$ be the corresponding group scheme of symplectic similitudes, whose
points in a $\Z$-algebra $A$ are given by
$G(A)=\GSp_{4,\Z}(A):=\{M\in \mathrm{M}_4(A)|\exists\alpha(M)\in A^*
\textrm{ with }MJM^{-1}=\alpha(M).J\}.$

Let $h_i^S:\Sb\to G_\R$ be the morphism that maps
$z=a+ib\in \C^*=\Sb(\R)$
to the matrix
$$
\left(\begin{array}{cc}
aI_2 & bI_2\\
-bI_2 & aI_2
\end{array}\right),
$$
and denote by $\Sc^\pm$ the $G(\R)$-conjugacy class of $h_i$,
that is usually called the two dimensional Siegel space.

\begin{definition}
The datum $(\GSp_4,\Sc^\pm)$ is called the \emph{Siegel Shimura datum}.
\end{definition}

Let $K$ be the compact open subgroup $G(\hat{\Z})$ of $G(\A_f)$.
The associated Shimura variety is by definition
$$
\begin{array}{rl}
\Sh_K(\GSp_4,\Sc^\pm)&
=G(\Q)\backslash (\Sc^\pm\times G(\A_f)/K)\\
& =\GSp_4(\Z)\backslash \Sc^\pm,
\end{array}
$$
i.e., the Siegel modular variety, which is the moduli space of principally
polarized Abelian surfaces.

Let $h_0^\Sc:\G_{m,\C}^2\to G_\C$ be the morphism that associates to
$(x,y)\in(\C^*)^2$ the matrix
$$
\left(\begin{array}{cccc}
a & 0 & b & 0\\
0 & x & 0 & 0\\
-b& 0 & a & 0\\
0 & 0 & 0 & y
\end{array}\right),
$$
with $a=\frac{x+y}{2}$ and $b=\frac{x-y}{2i}$.
Let $\Rc_\Sc:=G(\R)\cdot h_0^\Sc$ be the $G(\R)$-conjugacy class of
$h_0^\Sc$. The centralizer of $h_0^\Sc$ is a (non-split) maximal torus of 
$G(\R)$ (the subgroup $T(\R)$ of $\C^*\times (\R^*)^2$
given by $z\bar{z}=xy$), which implies that
$\Rc_\Sc\cong \GSp_4(\R)/T(\R)$.

The Shimura space
$\Sh_{G(\hat{\Z})}(G,\Rc_\Sc)$
is a moduli space for tuples
$$(V,F_x,F_y,\Pi,F,\psi)$$
where $V$ is a free $\Z$-module of rank $4$, $F_x$, $F_y$ are two distinct
real lines and $\Pi$ is a real plane in $V_\R$ such that
$V_\R=F_x\oplus F_y\oplus \Pi$. Moreover, $F$ is a complex line in $\Pi_\C$ such
that $\Pi_\C=F\oplus \overline{F}$, and $\psi:V\times V\to \Z$ is a symplectic
form that respects the decomposition
$$V_\C=(F_x\oplus F)\bigoplus (F_y\oplus\overline{F}).$$

This new kind of linear algebra object is quite strange and does not seem
to have an easy non-commutative geometric interpretation because it mixes usual
complex structures with foliations on tori. Such objects appear, however, quite
often in number theory.

\begin{example}
Let $K:=\Q[x]/(x^4-2)$. We have a decomposition
$K\otimes_\Q\R\cong \R_x\times \R_y \times \C_z$, and the $\R$-algebra
$\C_z$ decomposes over $\C$ as a product of two copies of $\C$. Fixing a
nice alternating form on $K$ gives
us exactly a point in $\Sh_{G(\hat{\Z})}(G,\Rc_\Sc)$. This point is
called a \emph{special point} or a \emph{quadratic multiplication point}.
We showed in \cite{expose-rivages} that the counting of points
of this type in the Shimura space involves interesting number theoretical
information, as in the case of geodesics studied in Section
\ref{special-geodesics-CFT}.
\end{example}

\section{Some open problems}
Here are some open problems\footnote{Some of them are collected from
the literature.}
in our work.
\begin{itemize}
\item Find a higher dimensional and/or algebraic analog of non-commutative 2-tori
adapted to number theoretical purposes.
\item Understand, in the case of geodesics, the relation of our
work with Manin's quantum theta functions, and Stark numbers.
\item More generally, find an adelic formulation of Stark's conjectures
for quadratic fields over totally real fields.
\item Define a well-behaved and not ad hoc notion of level structure
on non-commutative tori (resp., on Polishchuk's t-structures). 
\item Find the good higher dimensional analogs of Polishchuk's t-structures
on categories of coherent sheaves of Abelian varieties.
\item Study moduli spaces for stability conditions on these categories.
\item Clarify, if it exists, the relationship with Darmon's work on Stark's
conjecture for real quadratic fields.
\end{itemize}


\begin{thebibliography}{DMOS82}

\bibitem[AKO02]{Andre-Kahn}
Yves Andr{\'e}, Bruno Kahn, and Peter O'Sullivan.
\newblock Nilpotence, radicaux et structures mono\"\i dales.
\newblock {\em Rend. Sem. Mat. Univ. Padova}, 108:107--291, 2002.

\bibitem[BC95]{Bost-Connes}
J.-B. Bost and A.~Connes.
\newblock Hecke algebras, type {III} factors and phase transitions with
  spontaneous symmetry breaking in number theory.
\newblock {\em Selecta Math. (N.S.)}, 1(3):411--457, 1995.

\bibitem[BL01]{BoJi}
Armand Borel and Ji~Lizhen.
\newblock {\em Compactifications of locally symmetric spaces}.
\newblock IAS/Michigan, 2001.
\newblock Prepublication, june 2001.

\bibitem[Car98]{Cartier}
Pierre Cartier.
\newblock La folle journ\'ee, de {G}rothendieck \`a {C}onnes et {K}ontsevich:
  \'evolution des notions d'espace et de sym\'etrie.
\newblock In {\em Les relations entre les math\'ematiques et la physique
  th\'eorique}, pages 23--42. Inst. Hautes \'Etudes Sci., Bures, 1998.

\bibitem[CDS97]{Connes1}
Alain Connes, Michael.R. Douglas, and Albert Schwarz.
\newblock Noncommutative geometry and matrix theory: compactification on tori.
\newblock {\em arXiv}, (http://fr.arXiv.org/abs/hep-th/9711162), 1997.

\bibitem[CM95]{Connes5}
A.~Connes and H.~Moscovici.
\newblock The local index formula in noncommutative geometry.
\newblock {\em Geom. Funct. Anal.}, 5(2):174--243, 1995.

\bibitem[Con80]{Connes4}
Alain Connes.
\newblock {$C\sp{\ast} $} alg\`ebres et g\'eom\'etrie diff\'erentielle.
\newblock {\em C. R. Acad. Sci. Paris S\'er. A-B}, 290(13):A599--A604, 1980.

\bibitem[Con85]{Connes3}
Alain Connes.
\newblock Noncommutative differential geometry.
\newblock {\em Inst. Hautes \'Etudes Sci. Publ. Math.}, (62):257--360, 1985.

\bibitem[Con99]{Connes6}
Alain Connes.
\newblock Trace formula in noncommutative geometry and the zeros of the
  {R}iemann zeta function.
\newblock {\em Selecta Math. (N.S.)}, 5(1):29--106, 1999.

\bibitem[Del74]{De11}
Pierre Deligne.
\newblock Th\'eorie de {H}odge. {III}.
\newblock {\em Inst. Hautes \'Etudes Sci. Publ. Math.}, (44):5--77, 1974.

\bibitem[Del79]{De4}
Pierre Deligne.
\newblock Vari\'et\'es de {S}himura: interpr\'etation modulaire, et techniques
  de construction de mod\`eles canoniques.
\newblock In {\em Automorphic forms, representations and $L$-functions (Proc.
  Sympos. Pure Math., Oregon State Univ., Corvallis, Ore., 1977), Part 2},
  pages 247--289. Amer. Math. Soc., Providence, R.I., 1979.

\bibitem[DMOS82]{De1}
Pierre Deligne, James~S. Milne, Arthur Ogus, and Kuang-yen Shih.
\newblock {\em Hodge cycles, motives, and {S}himura varieties}.
\newblock Springer-Verlag, Berlin, 1982.
\newblock Philosophical Studies Series in Philosophy, 20.

\bibitem[GBVF01]{ENCG}
Jos{\'e}~M. Gracia-Bond{\'{\i}}a, Joseph~C. V{\'a}rilly, and H{\'e}ctor
  Figueroa.
\newblock {\em Elements of noncommutative geometry}.
\newblock Birkh\"auser Advanced Texts: Basler Lehrb\"ucher. [Birkh\"auser
  Advanced Texts: Basel Textbooks]. Birkh\"auser Boston Inc., Boston, MA, 2001.

\bibitem[Lod98]{Loday}
Jean-Louis Loday.
\newblock {\em Cyclic homology}, volume 301 of {\em Grundlehren der
  Mathematischen Wissenschaften [Fundamental Principles of Mathematical
  Sciences]}.
\newblock Springer-Verlag, Berlin, second edition, 1998.
\newblock Appendix E by Mar\'\i a O. Ronco, Chapter 13 by the author in
  collaboration with Teimuraz Pirashvili.

\bibitem[Man]{Manin3}
Yuri~I. Manin.
\newblock {Real Multiplication and noncommutative geometry}.

\bibitem[Mar04]{Mat4}
Matilde Marcolli.
\newblock {Lectures on Arithmetic Noncommutative Geometry}.
\newblock {\em arXiv}, (math.QA/0409520), 2004.

\bibitem[Mil83]{Milne5}
J.~S. Milne.
\newblock The action of an automorphism of {${\bf C}$} on a {S}himura variety
  and its special points.
\newblock In {\em Arithmetic and geometry, Vol. I}, volume~35 of {\em Progr.
  Math.}, pages 239--265. Birkh\"auser Boston, Boston, MA, 1983.

\bibitem[Mil90]{Mi3}
J.~S. Milne.
\newblock Canonical models of (mixed) {S}himura varieties and automorphic
  vector bundles.
\newblock In {\em Automorphic forms, Shimura varieties, and $L$-functions,
  Vol.\ I (Ann Arbor, MI, 1988)}, pages 283--414. Academic Press, Boston, MA,
  1990.

\bibitem[MM01]{Manin-Marcolli}
Yuri~I. Manin and Matilde Marcolli.
\newblock Holography principle and arithmetic of algebraic curves.
\newblock {\em Adv. Theor. Math. Phys.}, 5(3):617--650, 2001.

\bibitem[MM02]{Manin-Marcolli2}
Yuri~I. Manin and Matilde Marcolli.
\newblock Continued fractions, modular symbols, and noncommutative geometry.
\newblock {\em Selecta Math. (N.S.)}, 8(3):475--521, 2002.

\bibitem[Pau04]{expose-rivages}
F~Paugam.
\newblock Quelques bords irrationnels de vari\'et\'es de shimura.
\newblock {\em Universitaetsverlag Goettingen, Mathematisches Institut,
  Seminars}, pages 1--12, 2004.

\bibitem[Pin90]{Pink}
Richard Pink.
\newblock {\em Arithmetical compactification of mixed {S}himura varieties},
  volume 209 of {\em Bonner Mathematische Schriften [Bonn Mathematical
  Publications]}.
\newblock Universit\"at Bonn Mathematisches Institut, Bonn, 1990.
\newblock Dissertation, Rheinische Friedrich-Wilhelms-Universit\"at Bonn, Bonn,
  1989.

\bibitem[Pol02]{Poli1}
Alexander Polishchuk.
\newblock Noncommutative 2-tori with real multiplication as noncommutative
  projective varieties.
\newblock {\em arXiv}, (http://fr.arXiv.org/abs/math.AG/0212306), 2002.

\bibitem[Pol03]{Poli2}
Alexander Polishchuk.
\newblock Classification of holomorphic vector bundles on noncommutative
  two-tori.
\newblock {\em arXiv}, (http://fr.arXiv.org/abs/math.QA/0308136), 2003.

\bibitem[SGA64]{SGA3III}
{\em Sch\'emas en groupes. {I}{I}{I}: {S}tructure des sch\'emas en groupes
  r\'eductifs}.
\newblock Springer-Verlag, Berlin, 1962/1964.

\bibitem[SW99]{Seiberg-Witten}
Nathan Seiberg and Edward Witten.
\newblock String theory and noncommutative geometry.
\newblock {\em J. High Energy Phys.}, (9):Paper 32, 93 pp. (electronic), 1999.

\end{thebibliography}

\end{document}